\documentstyle[11pt]{article}  

\newcommand{\vm}{\vspace{0.2cm}}  
  
\begin{document}

\hsize=16.0truecm\vsize=22.5truecm\vglue6.3truecm


\begin{center}

 {\Large {\bf A Strategy for Proving Riemann Hypothesis}}

\vspace{.7cm}
M. Pitk\"anen$^1$\\ 
\vspace{.7cm}
{\footnotesize $^1$  Department of Physical Sciences, High Energy Physics Division,\\ PL 64, FIN-00014, University of Helsinki, Finland.\\ 
\indent matpitka@rock.helsinki.fi, http://www.physics.helsinki.fi/$\sim$matpitka/.\\ \indent Recent address: Kadermonkatu 16,10900, Hanko, Finland.}

\end{center}
\vspace{.2cm}

\noindent {\bf Abstract}.  
A strategy for proving  Riemann hypothesis is suggested.  The vanishing of the
Rieman Zeta reduces to an orthogonality condition for the eigenfunctions of a non-Hermitian operator $D^+$ having the zeros of Riemann Zeta as its eigenvalues. The construction of $D^+$ is inspired by the  conviction  that Riemann Zeta is associated with a physical system allowing conformal transformations as  its  symmetries. The eigenfunctions  
of $D^+$ are analogous to the so called  coherent states and in general not orthogonal  to each other. The states orthogonal to a vacuum state (which
has a negative norm squared) correspond  to the zeros of the Riemann Zeta. The 
induced metric in the space ${\cal{V}}$ of states which
correspond  to the zeros of the Riemann Zeta at the critical line $Re[s]=1/2$
is hermitian and both  hermiticity 
and positive definiteness properties imply  Riemann hypothesis.  
Conformal invariance in the sense of gauge invariance
allows only the states belonging to  ${\cal{V}}$. 
Riemann hypothesis follows  also  from a 
restricted form of a dynamical  conformal invariance in ${\cal{V}}$
and one can reduce the proof to a  standard analytic argument used in Lie
group theory.  

\newpage

\section{Introduction}
\label{intro}

The Riemann hypothesis \cite{Riemann,Titchmarch}
 states that the non-trivial  zeros (as opposed
to zeros at $s= -2n$, $n\geq 1$ integer) of Riemann Zeta function obtained by
analytically continuing the function 

\begin{eqnarray}
\zeta (s) = \sum_{n=1}^{\infty} \frac{1}{n^s} 
\end{eqnarray}

\noindent from the region $Re[s]>1$
to the  entire complex plane,  lie on the line  $Re[s]=1/2$. 
Hilbert and Polya  conjectured a long time ago
that the non-trivial zeroes of Riemann Zeta function could have  spectral
interpretation in terms of the eigenvalues of a suitable self-adjoint differential operator $H$ such that the eigenvalues of
this operator correspond to the imaginary
parts  of   the  nontrivial zeros $z=x+iy$
of $\zeta$. One can however consider a variant of this  hypothesis
stating that  the eigenvalue
spectrum of a   non-hermitian operator   $D^+$
contains  the non-trivial zeros of $\zeta$. 
The eigenstates in question are 
eigenstates of an annihilation operator type operator $D^+$
and analogous to the  so called coherent states encountered in 
quantum physics \cite{field}. In particular, the eigenfunctions 
are in general non-orthogonal and this is a quintessential element of the
the proposed strategy of proof.

In the following an explicit operator having as its
eigenvalues the non-trivial zeros of $\zeta$ is constructed.

a) The construction relies
crucially on the interpretation of the vanishing of
$\zeta$  as an orthogonality condition in a  hermitian metric which
is is a priori more general than Hilbert space inner product.

b) Second basic element is  the scaling invariance motivated by the belief
that $\zeta$ is associated with a physical
system which has superconformal transformations \cite{sconf} as its symmetries.

The core elements of the construction are following.

a)  All complex
numbers are  candidates for the eigenvalues of $D^+$ (formal hermitian
conjugate of $D$)
and   genuine
eigenvalues are selected by the requirement that the
condition  $D^{\dagger}=D^+$
holds true in the set of the genuine eigenfunctions. This condition
is equivalent with the hermiticity of the  metric defined
by a function proportional to $\zeta$.

b)  The  eigenvalues turn out to consist of $z=0$ and the
non-trivial zeros of $\zeta$ and only the
eigenfunctions  corresponding to the  zeros  with $Re[s]= 1/2$ define
a subspace possessing a hermitian metric.
The vanishing of $\zeta$  tells that the 'physical'  positive norm
eigenfunctions (in general {\it not} orthogonal to each other),   are
orthogonal to the   'unphysical'  negative norm  eigenfunction
associated with the eigenvalue $z=0$.

The proof of the Riemann hypothesis by reductio ad absurdum
results if one assumes that the space ${\cal{V}}$ spanned by the  states 
corresponding to the zeros of $\zeta$ inside the critical strip
has a hermitian induced metric.  Riemann hypothesis follows also from
the requirement that the induced metric in the spaces 
subspaces ${\cal{V}}_s$ of ${\cal{V}}$ spanned
by the states $\Psi_s$ and $\Psi_{1-\overline{s}}$ does not possess
negative eigenvalues: this condition is
equivalent with the positive definiteness of the metric
in ${\cal{V}}$.  Conformal invariance in the sense of gauge invariance
allows only the states belonging to  ${\cal{V}}$. 
Riemann hypothesis follows  also  from a 
restricted form of a dynamical  conformal invariance in ${\cal{V}}$.
This allows the reduction of the  proof to a  standard  analytic
argument used in Lie-group theory.

\section{Modified form of the Hilbert-Polya conjecture}
\label{sec1:1}

 One can  modify  the Hilbert-Polya  conjecture by assuming 
scaling invariance and giving up
the hermiticity of the Hilbert-Polya operator. This means
introduction of the    non-hermitian operators $D^+$ and $D$
which are hermitian conjugates of each other such that
 $D^{+}$ has the nontrivial zeros of $\zeta$
as its  complex eigenvalues

\begin{eqnarray}
\begin{array}{ll}
D^+ \Psi= z \Psi  .\\ 
\end{array}
\end{eqnarray}

\noindent The counterparts of the  so called  coherent states
\cite{field} are in question and 
the eigenfunctions  of $D^+$ are not expected to be orthogonal in general.
  The following construction
is based on the idea that
$D^{+}$  also allows  the eigenvalue $z=0$ and 
that the vanishing of $\zeta$ at $z$ expresses the  orthogonality
of the states with eigenvalue $z=x+iy\neq 0$ and the state
with eigenvalue  $z=0$ which turns out to have a negative norm.

  The trial 

\begin{eqnarray}
\begin{array}{ll}
D  = L_0+ V   ,  &D^+ = -L_0+V\\
\\
L_0=t\frac{d}{dt} ,& V= \frac{dlog(F)}{d(log(t))} =t\frac{dF}{dt}\frac{1}{F}  \\ 
\end{array}
\end{eqnarray}

\noindent is motivated by the requirement of   invariance
with respect to scalings $t\rightarrow \lambda t$ and
$F\rightarrow \lambda F$. The range of variation for
the variable $t$ consists of non-negative real numbers $t\geq 0$. 
The scaling invariance implying
conformal invariance  (Virasoro generator $L_0$ represents
scaling which plays a  fundamental role in the
superconformal theories \cite{sconf})
is motivated by the  belief that $\zeta$ codes for  the physics
of a quantum critical system having, not only supersymmetries \cite{Berry}, 
but  also superconformal transformations as its basic symmetries
(see the chapter "Riemann Hypothesis" of \cite{TGD}).

\section{Formal solution of the eigenvalue equation for operator $D^{+}$}
\label{sec2:1}
  One can formally solve the eigenvalue equation 

\begin{eqnarray}
D^+ \Psi_z= \left[-t\frac{d}{dt} +t\frac{dF}{dt}\frac{1}{F}\right]\Psi_z = 
 z \Psi_z  .
\end{eqnarray}

\noindent for  $D^{+}$ by  factoring the eigenfunction to a product: 

\begin{eqnarray}
\Psi_z= f_z F .
\end{eqnarray}

\noindent The substitution into the eigenvalue equation gives

\begin{eqnarray}
L_0 f_z= t\frac{d}{dt}f_z = -zf_z
\end{eqnarray}

\noindent allowing as its solution the functions

\begin{eqnarray}
f_z(t) = t^{z}  .
\end{eqnarray}

\noindent These functions are nothing but eigenfunctions of the scaling
operator $L_0$ of the superconformal
algebra analogous to the eigenstates of a translation operator.
 A priori all complex numbers
$z$ are candidates for the eigenvalues of $D^+$  and one must select
the genuine eigenvalues by applying  the requirement
$D^{\dagger}=D^+$ in the space spanned by the genuine eigenfunctions.

It must be emphasized  that $\Psi_z$ is {\it not} an  eigenfunction of
$D$. Indeed, one has

\begin{eqnarray}
D\Psi_z = -D^+ \Psi_z + 2V\Psi_z = z\Psi_z + 2V\Psi_z .
\end{eqnarray}

\noindent  This is in accordance with the   analogy
with the  coherent states
which are eigenstates of annihilation operator but not those
of creation operator.

\section{$D^+=D^{\dagger}$ condition and hermitian form}
\label{sec3:1}
The requirement that $D^{+}$ is indeed  the hermitian conjugate
of $D$ implies that the hermitian form
satisfies

\begin{eqnarray}
\langle f\vert D^{+} g\rangle   = \langle Df\vert  g \rangle  . 
\end{eqnarray}

\noindent This condition implies

\begin{eqnarray}
\langle \Psi_{z_1}\vert D^+\Psi_{z_2}\rangle
= \langle D\Psi_{z_1}\vert \Psi_{z_2}\rangle
 .
\end{eqnarray}

The first (not quite correct) guess is that the hermitian form 
is defined as an  integral of the product $\overline{\Psi}_{z_1}\Psi_{z_2}$
of the eigenfunctions of the operator $D$ over 
the non-negative real axis using a suitable
integration measure. 
 The hermitian form  can be defined by continuing
the integrand from the non-negative
real axis to the entire  complex $t$-plane  and noticing that it has a cut
along the non-negative real axis. This suggests  the definition
of the hermitian form, not as a mere integral over the non-negative real
axis, but as a contour integral along   curve $C$ defined so that
it encloses the non-negative real axis, that is $C$

a) traverses the non-negative real axis along the line $Im[t]=0_-$
from  $t=\infty + i0_-$  to $t=0_+ +i0_-$, 

b) encircles the origin around a small circle
from $t=0_+ +i0_-$   to $t= 0_+ +i0_+$,

c) traverses the  non-negative real axis
along the line $Im[t]=0_+$  from $t=0_+ +i0_+$ to 
$t= \infty+i0_+$ . \\
Here $0_{\pm}$  signifies taking the limit $x= \pm\epsilon$, $\epsilon>0$, $\epsilon \rightarrow 0$.  

 $C$ is the correct choice if   the integrand  defining the
inner product approaches zero sufficiently fast at
the limit  $Re[t]\rightarrow \infty$.  Otherwise
one must assume  that    the integration contour 
 continues along the  circle $S_R$ of radius $R\rightarrow \infty$ 
back to $t=\infty+i0_-$ to form a closed contour. It however turns
out that this is not necessary. 
One can deform the integration contour 
rather  freely: the only constraint is 
 that the deformed integration contour  does not cross over any cut or pole
associated with the analytic continuation of   the integrand
 from the non-negative real axis to the entire
complex plane.

Scaling invariance  dictates the form of the  integration
measure appearing in the hermitian form  uniquely to be $dt/t$.
The hermitian form  thus obtained also makes  possible
to satisfy the crucial $D^+ = D^{\dagger}$
condition. The hermitian form is thus defined as

\begin{eqnarray}
\langle \Psi_{z_1}\vert \Psi_{z_2}\rangle   =  -\frac{K(z_{12})}{2\pi i}
 \int_C \overline{\Psi_{z_1}}\Psi_{z_2} \frac{dt}{t} .
\end{eqnarray}

\noindent        $K(z_{12})$ is real from the hermiticity
requirement and the behaviour as a function of
$z_{12}=z_1+\overline{z}_2$ by the requirement that the resulting
Hermitian form defines a positive definite inner product. The
value of  $K(1)$ can can be fixed by requiring that the states
corresponding to the zeros of $\zeta$ at the critical line have
unit norm: with this choice the vacuum state corresponding to
$z=0$ has negative norm. Physical intuition suggests that $K(z_{12})$ is
responsible for the Gaussian overlaps of the coherent states and
this suggests the behaviour

\begin{eqnarray}
K(z_{12})= exp(-\alpha \vert z_{12}\vert^2) ,
\end{eqnarray}

\noindent for which overlaps between states at critical line are
 proportional to $exp(-\alpha (y_1-y_2)^2)$ so that for $\alpha>0$ Schwartz  inequalities are certainly satisfied 
for large values of $\vert y_{12}\vert$.
Small values of $y_{12}$
are dangerous in this respect but since the matrix elements of the  metric
decrease   for small values of $y_{12}$ even for $K(z_{12})=1$,
it is possible to satisfy Schwartz inequalities for sufficiently large
value of $\alpha$.  It must be emphasized that the 
detailed behaviour
of $K$ is {\it not} crucial for the arguments relating to Riemann
hypothesis.

The possibility to
deform the shape of $C$ in wide limits realizes conformal invariance
stating that the change of the shape of the integration contour
induced by a conformal transformation, which is  nonsingular 
inside the integration  contour, leaves the
value of  the contour integral of an analytic function unchanged.
This scaling invariant hermitian form
is indeed a correct guess. By applying  partial integration one can write

\begin{eqnarray}
\langle \Psi_{z_1}\vert D^+\Psi_{z_2}\rangle
= \langle D\Psi_{z_1}\vert \Psi_{z_2}\rangle - \frac{K(z_{12})}{2\pi i}
\int_C dt \frac{d}{dt} \left[\overline{\Psi}_{z_1}(t) \Psi_{z_2}(t)\right]
 .
\end{eqnarray}

\noindent  The integral of a total differential comes from the operator 
$L_0=td/dt$ and must vanish.  For a  non-closed  integration contour $C$
the  boundary terms from the partial integration could spoil the  $D^+=D^{\dagger}$ condition unless the eigenfunctions vanish at the end points of the integration contour ($t=\infty + i0_{\pm}$).   

The explicit expression of the hermitian form is  given by

\begin{eqnarray}
\langle \Psi_{z_1}\vert \Psi_{z_2}\rangle
&=&  -\frac{K(z_{12})}{2\pi i}  \int_{C}  \frac{dt}{t} F^2(t) t^{z_{12}}    ,\nonumber\\
z_{12}&=&\overline{z}_1+z_2 .
\label{quadr} 
\end{eqnarray}

\noindent It must be emphasized that it is  
$\overline{\Psi}_{z_1}\Psi_{z_2}$ rather than
eigenfunctions which is continued from the 
non-negative real axis to the  complex $t$-plane: therefore one
indeed obtains an analytic function as a result.
 
 An essential role in the argument claimed to prove the
Riemann hypothesis  is played by 
the crossing symmetry

\begin{eqnarray}
\langle\Psi_{z_1}\vert \Psi_{z_2} \rangle  & =&  
  \langle\Psi_0\vert \Psi_{\overline{z}_1+z_2} \rangle  
\label{cross}
\end{eqnarray}

\noindent of the hermitian form.  This symmetry is analogous to 
the crossing symmetry of particle physics  stating that 
 the S-matrix is symmetric with respect to the  
replacement  of the particles in the  initial state with
their antiparticles in  the final state or vice versa \cite{field}.

The hermiticity of the hermitian form   implies

\begin{eqnarray}
\langle \Psi_{z_1}\vert \Psi_{z_2}\rangle
&=&  
\overline{\langle \Psi_{z_2}\vert \Psi_{z_1}\rangle} .
\label{herm} 
\end{eqnarray}

\noindent This condition, which is {\it not} trivially satisfied,
  in fact determines the eigenvalue spectrum.

\section{How to choose the function $F$?}
\label{sec4:1}

 The remaining task is to choose the function $F$ in such a manner
that the orthogonality conditions for the solutions $\Psi_0$ and $\Psi_z$
reduce to the condition that $\zeta$  or some
function proportional to $\zeta$ vanishes at the point $-z$. 
The definition of   $\zeta$   based
on   analytical continuation performed 
by Riemann suggests how to proceed.
Recall that the expression of  $\zeta$  converging in
the region $Re[s]>1$ following from the basic definition
of $\zeta$ and elementary properties of $\Gamma$ function
\cite{Titchmarch} reads   as

\begin{eqnarray}
\Gamma (s) \zeta (s)= \int_{0}^{\infty} \frac{dt}{t} \frac{exp(-t)}{\left[1-exp(-t)\right]} t^s 
 . 
\end{eqnarray}
 
\noindent One can analytically continue this expression to
a function defined  in  the entire complex plane by 
noticing that the  integrand is discontinuous along the  cut 
extending from $t=0$ to $t=\infty$. 
Following Riemann it is however more convenient to consider the
discontinuity for a function obtained by multiplying the integrand
with  the factor

 $$(-1)^s \equiv  exp(-i\pi s)   .$$

 The discontinuity  $Disc (f)\equiv f(t)-f(texp(i2\pi))$  of the
resulting function is given by

\begin{eqnarray}
Disc\left[\frac{exp(-t)}{\left[1-exp(-t)\right]} (-t)^{s-1}\right]= 
-2i sin(\pi s)  \frac{exp(-t)}{\left[1-exp(-t)\right]} t^{s-1}  .
\end{eqnarray}

\noindent  The discontinuity 
 vanishes  at the limit $t\rightarrow 0$ for $Re[s]>1$.
 Hence one can  define  $\zeta$ by modifying the integration
contour from the non-negative  real axis to an integration
contour $C$ enclosing non-negative real axis defined in the
previous section.

This amounts  to writing 
the analytical continuation of $\zeta (s)$ in the form  

\begin{eqnarray}
-2i\Gamma (s) \zeta (s) sin(\pi s)= 
\int_{C} \frac{dt}{t} \frac{exp(-t)}{\left[1-exp(-t)\right]} 
(-t)^{s-1}  . 
\label{riemann}
\end{eqnarray}

\noindent This expression equals to $\zeta (s)$  for $Re[s]>1$
and defines $\zeta(s)$ in the entire complex plane since 
the integral around the origin  eliminates the singularity.

The crucial observation is that the 
integrand on the righthand side of Eq.  \ref{riemann}
has precisely the same general form as that appearing
in  the hermitian form defined  in Eq. \ref{quadr} defined 
using the same integration contour $C$. 
The integration measure is $dt/t$, 
the factor $t^{s}$ is  of the same form as 
the factor $t^{\overline{z}_1 +z_2 }$ appearing 
in the hermitian form,  and the function $F^2(t)$ is given by

 $$F^2(t)= \frac{exp(-t) }  {1-exp(-t) }  .$$

\noindent Therefore one can make the identification

\begin{eqnarray}
F(t) =  \left[ \frac{exp(-t)}{1-exp(-t)}\right]^{1/2}  . 
\end{eqnarray}

\noindent Note that the argument of
the  square root is non-negative on the non-negative
real axis and that $F(t)$
decays exponentially on the non-negative real axis and
has $1/\sqrt{t}$ type singularity at origin. From this it
follows that the eigenfunctions $\Psi_z (t)$ approach zero exponentially
at the limit $Re[t]\rightarrow \infty$ so that one can use the
non-closed integration contour $C$.

With this assumption,  the hermitian form reduces to the  expression

\begin{eqnarray}
\langle\Psi_{z_1}\vert \Psi_{z_2} \rangle   &=& -\frac{K(z_{12})}{2\pi i } 
\int_{C} \frac{dt}{t} \frac{exp(-t)}{\left[1-exp(-t\right]} 
(-t)^{z_{12}}\nonumber\\
\nonumber\\
&=& \frac{K(z_{12})}{\pi }
 sin(\pi z_{12})   \Gamma (z_{12})  \zeta (z_{12})  .   
\label{inner}
\end{eqnarray}

\noindent Recall that the definition $z_{12}= \overline{z}_1+z_2$ is
adopted.    Thus the orthogonality of the eigenfunctions is equivalent
to the vanishing of  $\zeta (z_{12})$ if $K(z_{12})$ is positive definite.

 \section{Study of the hermiticity condition}
\label{sec5:1}

In order to derive information about the
spectrum one must explicitely study what the statement
that $D^{\dagger}$ is hermitian conjugate of $D$ means.
The defining equation is just the generalization
of the  equation

\begin{eqnarray}
A^{\dagger}_{mn}= \overline{A}_{nm} .
\end{eqnarray}

\noindent defining the notion of hermiticity for matrices. 
Now indices $m$ and $n$ correspond to the eigenfunctions
$\Psi_{z_i}$,  and one obtains

$$\begin{array}{l}
\langle \Psi_{z_1}\vert D^+\Psi_{z_2}\rangle=z_2\langle \Psi_{z_1}\vert \Psi_{z_2}\rangle
= \overline{\langle \Psi_{z_2}\vert D\Psi_{z_1}\rangle }
=\overline{\langle D^+\Psi_{z_2}\vert \Psi_{z_1}\rangle}
= z_2\overline{\langle \Psi_{z_2}\vert \Psi_{z_1}\rangle} .\\
\end{array}$$

\noindent   Thus  one has

\begin{eqnarray}
G(z_{12}) &=& \overline{G(z_{21})}= \overline{G(\overline{z}_{12})}\nonumber\\
G(z_{12}) &\equiv& \langle \Psi_{z_1}\vert \Psi_{z_2}\rangle
 . 
\end{eqnarray}

\noindent  The condition  states that the hermitian form  defined
by the contour integral is
indeed hermitian. This is {\it not} trivially true.  
Hermiticity condition obviously determines the spectrum of 
the eigenvalues of $D^+$.

To see the implications of the  hermiticity condition,
one must study the behaviour of the function $G(z_{12})$ 
under complex conjugation
of both the argument and the value of the function itself. To achieve this
one must write the integral 

$$G(z_{12})   = -\frac{K(z_{12})}{2\pi i } 
\int_{C} \frac{dt}{t} \frac{exp(-t)}{\left[1-exp(-t)\right]} (-t)^{z_{12}}$$

\noindent  in a  form
from which one can easily deduce the behaviour of this function
under complex conjugation. To achieve this, one must
perform the change $t\rightarrow u=log(exp(-i\pi )t)$ of the integration variable
giving

\begin{eqnarray}
G(z_{12})   &=& -\frac{K(z_{12})}{2\pi i}  
\int_{D} du \frac{exp(-exp(u))  }{\left[1-exp(-(exp(u)))\right]} 
exp(z_{12}u)  . \nonumber\\
\end{eqnarray}

\noindent  Here $D$ denotes the image of the integration contour
$C$ under $t \rightarrow u= log(-t)$.  $D$
is a fork-like  contour which\\
a)  traverses the line  $Im[u]=i\pi$ from $u=\infty+i\pi$ to $u=-\infty +i\pi$ , \\
 b) continues from $-\infty +i\pi$ to $-\infty -i\pi$ along the imaginary
$u$-axis (it is easy to see that the contribution from this part
of the contour vanishes),\\
c) traverses the real $u$-axis from $u=-\infty-i\pi$ to $u=\infty-i\pi$,\\

The integrand differs on the line $Im[u]=\pm i\pi$ from
that on the line $Im[u]=0$ by the factor  $ exp(\mp i\pi z_{12})$ so that
one can write $G(z_{12})$ as integral over real $u$-axis

\begin{eqnarray}
G(z_{12})   &=& -\frac{K(z_{12})}{\pi }  
sin(\pi z_{12}) \int_{-\infty}^{\infty} du \frac{exp(-exp(u))  }{\left[1-exp(-(exp(u)))\right]} 
exp(z_{12}u)  . \nonumber\\
\end{eqnarray}

\noindent 
From this form the effect of the
transformation $G(z)\rightarrow \overline{G(\overline{z})}$
can be deduced. Since the integral is along the real $u$-axis, complex
conjugation amounts only to the replacement  $z_{21}\rightarrow z_{12}$,
and one has

\begin{eqnarray}
  \overline{G(\overline{z}_{12})}    &=& -\frac{\overline{K(z_{21})}}{\pi } 
 \times \overline{sin(\pi z_{21})} 
\int_{-\infty}^{\infty} du \frac{exp(-exp(u))  }{\left[1-exp(-(exp(u)))\right]} 
exp(z_{12}u)   \nonumber\\
&=& \frac{\overline{K(z_{21})}}{K(z_{12})}\times \frac{\overline{sin(\pi z_{21})}}{sin(\pi z_{12})}  G(z_{12})     .
\end{eqnarray}

\noindent Thus the hermiticity  condition reduces to the
condition

\begin{eqnarray}
G(z_{12}) = \frac{\overline{K(z_{21})}}{K(z_{12})}\times \frac{\overline{sin(\pi z_{21})}}{sin(\pi z_{12})} 
\times G(z_{12})   .
\end{eqnarray}

\noindent The reality of  $K(z_{12})$  
guarantees  that the  diagonal matrix elements  of the metric are  real.

For non-diagonal matrix elements 
there are two manners to satisfy the hermiticity condition.

a) The condition 
 
\begin{eqnarray}
G(z_{12}) = 0  
\end{eqnarray}

\noindent is the only  manner to satisfy the  hermiticity
condition for $x_1+x_2\neq n$, $y_1-y_2 \neq 0$.  This implies 
the vanishing of $\zeta$:

\begin{eqnarray}
\zeta (z_{12}) = 0  \mbox{~for~} 0<x_1+x_2<1    .
\end{eqnarray}

\noindent In particular, this condition
must be  true for $z_1=0$ and $z_2=1/2+iy$.
Hence  the physical states with
the  eigenvalue $z=1/2+iy$ must correspond to the zeros
of $\zeta$.

b) For the  non-diagonal matrix elements of the metric
the condition

\begin{eqnarray}
 exp(i\pi  (x_1+x_2)) = \pm 1 
\end{eqnarray}

\noindent guarantees the reality of $sin (\pi z_{12})$ factors. 
This requires 

\begin{eqnarray}
x_1+x_2=n  .
\end{eqnarray}

\noindent 
The highly non-trivial implication is that the the
vacuum state $\Psi_0$ and the zeros of $\zeta$ at the critical
line span a space having a hermitian.
Note that for $x_1=x_2=n/2$, $n\neq 1$, the diagonal matrix
elements of the metric vanish.

c) The metric is positive definite only
if the function $K(z_{12})$ decays sufficiently fast: this is due
to the exponential increase of the moduli of the matrix elements
$G(1/2+iy_1,1/2+iy_2)$  
for $K(z_{12})=1$ and for large values of $ \vert y_1-y_2\vert$ 
(basically due to the $sinh\left[\pi \left(y_1-y_2\right)\right]$-factor in
the  metric) 
implying the failure of the Schwartz inequality for
$\vert y_1-y_2\vert \rightarrow \infty$. 
Unitarity, guaranteing
probability interpretation  in quantum theory,  thus
requires that the parameter $\alpha$ characterizing the Gaussian decay of
$K(z_{12})=exp(-\alpha\vert z_{12}\vert^2)$ is  above some 
minimum value.

\section{Various assumptions implying  Riemann hypothesis}
\label{sec6:1}

As found, the general strategy for proving the Riemann hypothesis,
originally inspired by superconformal invariance,
leads to the construction of a set of eigenstates
for an operator $D^+$,  which is effectively an annihilation operator
acting in the space of complex-valued functions defined on the real half-line.
Physically the states are analogous to coherent
states and are not orthogonal to each other. The quantization of the
eigenvalues for the operator $D^+$
follows from the requirement that the metric, which is
defined by the integral defining the analytical
continuation of  $\zeta$,  and  thus  proportional to
$\zeta$ ($\langle s_1,s_2\rangle \propto \zeta (\overline{s}_1+s_2)$),
  is hermitian in   the  space
of the physical states.

 The nontrivial zeros of $\zeta$ are  known to  
belong to the critical strip defined by  $0<Re[s]<1$.  Indeed,
the theorem of Hadamard and de la Vallee  Poussin \cite{Edwards}
states the non-vanishing of $\zeta$ on the line $Re[s]=1$.
If $s$ is  a  zero of $\zeta$ inside the  critical strip,
then also $1-\overline{s}$ as well as $\overline{s}$ and $1-s$  are zeros.
If Hilbert space inner  product property
is not required so that the eigenvalues of the metric
tensor can be also negative in this subspace.  There
could be also unphysical zeros of $\zeta$ outside the critical
line $Re[s]=1/2$ but inside the critical strip
$0< Re[s]<1$. The problem is to
find whether the  zeros outside the critical line
are excluded, not only by the hermiticity but
also  by  the positive definiteness of the metric
necessary for the physical interpretation, and  perhaps also by
conformal invariance posed in some sense as a dynamical symmetry.
This turns out to be the case.

Before continuing it is convenient to introduce some notations.
Denote by  $\cal{V}$  the subspace
spanned by  $\Psi_s$ corresponding to 
the  zeros $s$ of $\zeta$ inside the critical strip, by ${\cal{V}}_{crit}$  the subspace  corresponding to  the  zeros  of $\zeta$ 
at the critical strip,  and by ${\cal{V}}_s$ the space spanned by the 
states $\Psi_s$ and $\Psi_{1-\overline{s}}$.
The basic idea behind the  following proposals is
that the basic objects of study are the spaces ${\cal{V}}$,
${\cal{V}}_{crit}$ and ${\cal{V}}_{s}$.

\subsection{How to restrict the metric to  ${\cal{V}}$?}
\label{sec61:2}

 One should somehow restrict the   metric defined in the space 
spanned by the states $\Psi_s$ labelled by  a continuous 
complex eigenvalue $s$  to the space
${\cal{V}}$ inside the critical strip spanned by a basis labelled
by discrete eigenvalues. Very naively, 
one could try to do this by simply putting
all other components of the metric to zero so that the states outside 
${\cal{V}}$ correspond to gauge degrees of freedom. This is consistent
with the interpretation of 
${\cal{V}}$ as a coset space formed by identifying states which differ
from each other by the addition of a superposition of states which
do not correspond to zeros of $\zeta$. 

An more elegant manner to realize the restriction of the metric
to ${\cal{V}}$  is to Fourier
expand states in the basis labelled by a complex number $s$
and  define the  metric in ${\cal{V}}$ using
double Fourier integral over the complex plane and
Dirac delta function restricting
the labels of both states to the set of zeros inside the critical strip:

\begin{eqnarray}
\langle \Psi^{1)}\vert\Psi^{2)}\rangle&= &\int d\mu (s_1)
 \int d\mu(s_2)\overline{\Psi}^{1)}_{s_1} \Psi^{2)}_{s_2} G (s_2+\overline{s}_1)  \delta (\zeta(s_1))  \delta (\zeta (s_2))\nonumber\\
&=& \sum_{\zeta (s_1)=0,\zeta (s_2)=0} \overline{\Psi}^{1)}_{s_1} \Psi^{2)}_{s_2}
G(s_2+\overline{s}_1)
\frac{1}{\sqrt{det(s_2)det(\overline{s}_1)}}   ,\nonumber\\
\nonumber\\
d\mu(s)&=& dsd\overline{s}  ,\hspace{1.0cm} 
 det(s)=\frac{ \partial (Re\left[\zeta (s)\right],Im\left[\zeta (s)\right])}{\partial (Re\left[s\right],Im\left[s\right])}.
\end{eqnarray}

\noindent Here the integrations are over the critical strip.
$det(s)$ is the Jacobian for the map $s\rightarrow \zeta (s)$
at $s$.  The appearence of the determinants might be crucial for
the  absence of negative norm states. The result means
that the metric $G_{{\cal{V}}}$ in ${\cal{V}}$
effectively reduces to a product

\begin{eqnarray}
G_{{\cal{V}}}&=& \overline{D} GD , \nonumber\\
D(s_i,s_j)&=& D(s_i)\delta (s_i,s_j), \nonumber\\
\overline{D}(s_i,s_j)&=& D(\overline{s_i})\delta (s_i,s_j)\nonumber\\
D(s) &=& \frac{1}{\sqrt{det(s)}} .
\end{eqnarray}

\noindent In the sequel the metric $G$  will be called reduced metric whereas 
$G_{{\cal{V}}}$ will be called the full metric. In fact, the
symmetry $D(s)= D(\overline{s})$ holds true
by the basic symmetries of $\zeta$ so that one has $D=\overline{D}$
and $G_{{\cal{V}}}= DGD$.   This  means   that Schwartz inequalities for
the eigen states of $D^+$ are not affected in the
replacement of $G_{{\cal{V}}}$ with  $G$. The two metrics can be
in fact transformed to each other by a mere scaling of the eigenstates
and are in this sense equivalent.

\subsection{Riemann hypothesis from the  hermicity of 
the metric in ${\cal{V}}$}
\label{sec62:2}

The mere requirement that the metric is hermitian in 
${\cal{V}}$  implies the Riemann hypothesis.  This can be
seen in the simplest manner as follows.
Besides the zeros   at the critical 
line $Re[s]=1/2$ also the symmetrically related
zeros inside critical strip have positive
norm squared but they do not have hermitian inner products
with the states at the critical line unless one
assumes that the inner product vanishes.
The assumption that the inner products between the states
at critical line and outside it vanish, implies additional
zeros of $\zeta$ and, by repeating the argument again and again,
one can fill the entire critical interval
$(0,1)$ with the zeros of $\zeta$ so that a reductio ad absurdum proof
for the Riemann hypothesis results.
 Thus the metric gives for   the 
states corresponding to the zeros of the
Riemann Zeta at the critical line  a special status as 
what might be called physical states. 

It should be noticed that the states in ${\cal{V}}_{s}$ and 
${\cal{V}}_{\overline{s}}$ have 
non-hermitian  inner products for $Re[s]\neq 1/2$ unless
these inner products vanish: for $Re[s]>1/2$ this however implies 
that $\zeta$ has a zero for $Re[s]>1$.

\subsection{Riemann hypothesis from the requirement
that the  metric in ${\cal{V}}$ is positive definite}

\label{sec63:2}

With a suitable choice of $K(z_{12})$
the metric is positive definite between states having  $y_1\neq y_2$. For 
$s$ and $1-\overline{s}$ one has $y_1=y_2$ implying  $K(z_{12})=1$
in ${\cal{V}}_s$. Thus the positive definiteness of the metric in ${\cal{V}}$ 
reduces to that for  the induced metric  in the spaces
${\cal{V}}_s$. This requirement  implies also Riemann hypothesis
as following argument shows.

The explicit expression for the norm  of a
$Re[s]=1/2$ state with respect to the full metric $G^{ind}_{\cal{V}}$
reads as

\begin{eqnarray}
G^{ind}_{{\cal{V}}}(1/2+iy_n, 1/2+iy_n)  
&=& D^2(1/2+iy) G^{ind}(1/2+iy_n, 1/2+iy_n), \nonumber\\
G^{ind}(1/2+iy_n, 1/2+iy_n)&=&  -\frac{K(z_{12})}{\pi }
 sin(\pi)  \Gamma (1) \zeta (1)  .
\end{eqnarray}

\noindent  Here $G^{ind}$   is the metric  in ${\cal{V}}_s$
induced from the  reduced metric $G$.
This expression involves formally a product of vanishing and infinite
factors and the value of expression must be defined as a limit by taking
in $Im[z_{12}]$ to zero. The requirement that the norm 
squared defined by $G^{ind}$ equals to one
 fixes the value of  $K(1)$:

\begin{eqnarray}
 K(1)&=& -\frac{\pi}{ sin(\pi)    \zeta (1)}=1  .
\end{eqnarray}

The components $G^{ind}$  in ${\cal{V}}_s$ are given by

\begin{eqnarray}
G^{ind}(s,s) &=& -
\frac{sin(2\pi Re[s])\Gamma (2Re[s])\zeta (2Re[s])}{ \pi} ,\nonumber\\
G^{ind}(1-\overline{s},1-\overline{s}) &=& -
\frac{sin(2\pi (1-Re[s]))\Gamma (2-2Re[s])\zeta (2(1-[Re[s]))}{\pi} ,\nonumber\\
G^{ind}(s,1-\overline{s})&=& G^{ind}(1-\overline{s},s)=  1 .
\end{eqnarray}

\noindent  
The determinant of the metric $G^{ind}_{{\cal{V}}}$
induced  from the  full metric 
reduces to the product 

\begin{eqnarray}
Det(G^{ind}_{{\cal{V}}})&=&D^2(s))D^2(1-\overline{s})\times Det(G^{ind})   .
\end{eqnarray}

\noindent Since the first factor is positive
definite,   it suffices to study the determinant of $G^{ind}$.
 At the limit $Re[s]=1/2$ 
 $G^{ind}$ formally reduces to

$$
\left(\begin{array}{rr} 
1 &1\\ 
1 &1\\
\end{array} \right) .
$$

\noindent This reflects the fact that the 
states $\Psi_s$ and $\Psi_{1-\overline{s}}$
are identical. The actual metric is of course
positive definite.  For  $Re[s]=0$ the 
$G^{ind}$ is  of the form 

$$
\left(\begin{array}{rr} 
-1 &1\\ 
1 &0\\
\end{array} \right) .
$$

\noindent The determinant
of $G^{ind}$ is negative so
that   the eigenvalues 
of both the full metric and reduced metric
are of opposite sign. 
The eigenvalues  for $G^{ind}$  are given  by
$(-1\pm \sqrt{5})/2$.

The determinant of $G^{ind}$  in ${\cal{V}}_s$ as a function of
$Re[s]$  is symmetric with respect to  $Re[s]=1/2$, equals to $-1$
at the end  points $Re[s]=0$ and $Re[s]=1$, 
and  vanishes at $Re[s]=1/2$. Numerical calculation shows that
the sign of the determinant of $G^{ind}$ 
  inside the interval $(0,1)$ is negative 
for $Re[s]\neq 1/2$. 
Thus the diagonalized form of the induced metric has
the  signature $(1,-1)$
except at the limit $Re[s]=1/2$,  when the  signature formally reduces to $(1,0)$. 
 Thus Riemann hypothesis follows if one can show   that
the metric induced to ${\cal{V}}_s$ does not allow  physical
states with a negative
norm squared. This requirement is  physically very natural. 
In fact, when the factor $K(z_{12})$ represents sufficiently 
rapidly vanishing Gaussian, this guarantees  the metric to ${\cal{V}}_{crit}$ has only non-negative eigenvalues.  Hence the positive-definiteness of
the metric, natural if there is real quantum system
behind the model, implies Riemann hypothesis.

\subsection{Riemann hypothesis and conformal invariance}

\label{sec64:2}

The basic  strategy for proving Riemann hypothesis has been based on 
the attempt to reduce Riemann hypothesis to 
invariance under  conformal  algebra or some
subalgebra of the conformal algebra in ${\cal{V}}$ 
or ${\cal{V}}_s$. That this kind of algebra
should act as a  gauge symmetry associated
with $\zeta$ is  very natural 
idea since conformal invariance is in
a well-defined sense the basic  symmetry group of complex analysis.

Consider now one particular  strategy based on conformal invariance
in the space of the eigenstates of $D^+$. 

{\it 1. Realization of conformal algebra as a spectrum generating algebra}

  The  conformal generators are realized
as operators  

\begin{eqnarray}
L_z &=& t^z D^+
\end{eqnarray}

\noindent act in the eigenspace of $D^+$
and obey the standard conformal  algebra without central
extension \cite{sconf}.
$D^+$ itself corresponds to the conformal generator $L_0$ acting as a scaling.
 Conformal generators obviously act as dynamical
symmetries transforming eigenstates of $D^+$  to each other.
 What is new is that now conformal weights $z$ have all 
possible complex values unlike in the standard case in which only integer
values are possible.  The vacuum state $\Psi_0$ having negative norm squared
is annihilated by the  conformal algebra so that the states orthogonal
to it (non-trivial zeros of $\zeta$ inside the critical
strip) form naturally another subspace which
should be conformally invariant in some sense. Conformal
algebra could act as gauge algebra
and some subalgebra of the conformal algebra  
 could act as a  dynamical symmetry.

{\it 2. Realization of conformal algebra as gauge symmetries}

 The definition of the metric in ${\cal{V}}$ involves in
an essential manner the mapping $s\rightarrow \zeta (s)$. This  suggests
that one should  define the 
gauge action of the conformal algebra
as

\begin{eqnarray}
\Psi_s &\rightarrow& \Psi_{\zeta (s)}\rightarrow 
L_z\Psi_{\zeta (s)}= \zeta_s\Psi_{\zeta (s)+z}\nonumber\\
&\rightarrow& 
\zeta_s\Psi_{\zeta^{-1}(\zeta (s)+z)} .
\end{eqnarray}

\noindent  Clearly, the action
involves a  map of the conformal weight $s$  to  $\zeta (s)$, the action
of  the conformal algebra to  $\zeta (s)$, 
and the mapping  of the transformed conformal weight $z+\zeta (s)$ back to
the complex plane by the  inverse of $\zeta$.  
The inverse image is in general non-unique 
but in case of ${\cal{V}}$ this does not matter since 
the  action annihilates automatically all  states in ${\cal{V}}$.
Thus conformal algebra indeed acts as a gauge symmetry. 
This symmetry does not however force Riemann hypothesis.

{\it 3. Realization of conformal algebra as dynamical symmetries}

 One can also study the action of the conformal
algebra  or its suitable sub-algebra
in ${\cal{V}}_s$ as a dynamical (as opposed
to gauge) symmetry realized as 

\begin{eqnarray}
\Psi_s&\rightarrow &L_z\Psi_s= s\Psi_{s+z}.
\end{eqnarray}
  
\noindent The states $\Psi_s$ and $\Psi_{1-\overline{s}}$ in ${\cal{V}}_s$
have  nonvanishing norms  and are obtained from each other by
the conformal generators $L_{1-2Re[s]}$ and  $L_{2Re[s]-1}$. 
For $Re[s]\neq 1/2$  the generators $L_{1-2Re[s]}$, $L_{2Re[s]-1}$, and $L_0$
generate $SL(2,R)$ algebra which is non-compact and generates
infinite number of states from the states of ${\cal{V}}_s$.
At the critical line this algebra reduces to the abelian
algebra spanned by  $L_0$.
The requirement that the  algebra naturally associated with
${\cal{V}}_s$ is a dynamical
symmetry and thus  generates only zeros of $\zeta$ leads to 
the conlusion that all points $s+n(1-2Re[s])$,
$n$ integer, must be zeros of $\zeta$. Clearly,
 $Re[s]=1/2$ is the only possibility
so that Riemann hypothesis follows. In this case the dynamical symmetry
indeed reduces to a gauge symmetry.

 There is clearly  a connection with the argument based on 
the requirement that the induced metric in ${\cal{V}}_s$  does not
possess negative eigenvalues. Since 
$SL(2,R)$  algebra acts as the isometries of the induced metric
for the zeros having  $Re[s]\neq 1/2$,
the  signature of the induced metric  must be $(1,-1)$.

{\it 4. Riemann hypothesis from the requirement that infinitesimal
isometries exponentiate}

 One could even try to prove  that
the entire subalgebra of the conformal algebra spanned by the generators
with conformal weights
$n(1-2Re[s])$ acts as a symmetry generating new zeros of $\zeta$
so that corresponding states are annihilated by gauge conformal algebra. 
If this holds, $Re[s]=1/2$ is the only possibility
so that Riemann hypothesis follows. In this case the dynamical
conformal symmetry
indeed reduces to a gauge symmetry.  

Since $L_{1-2Re[s]}$ acts as an infinitesimal 
isometry leaving the matrix element
$\langle \Psi_0\vert \Psi_s\rangle=0$ invariant, one can 
in spirit of Lie group theory argue that also
the exponentiated transformations 
$exp(tL_{1-2Re[s]})$   have the same property for all values of $t$. 
The exponential
action leaves $\Psi_0$ invariant and generates from $\Psi_s$ a superposition
of states with conformal weights $s+ n(1-2Re[s])$,  which all must be orthogonal
to $\Psi_0$ since $t$ is arbitrary. Since all zeros are inside the
critical strip,  $Re[s]=1/2$
is the only possibility.

A more explicit formulation of this idea is based on a first order differential
equation for the integral representation of $\zeta$. One can write  
the matrix element of the metric
 using the analytical continuation of $\zeta (s)$: 

\begin{eqnarray}
G(s) &=&-2i\Gamma (s) \zeta (s) sin(\pi s)= H(s,a)_{\vert a=0} , \nonumber\\
H(s,a) &=&\int_{C} \frac{dt}{t} \frac{exp(-t +a(-t)^{1-2x})}{\left[1-exp(-t)\right]} 
(-t)^{x+iy-1} . 
\label{riemann1}
\end{eqnarray}

\noindent If $s=x+iy$ is zero of $\zeta$ then also $1-x+iy$ is zero of $\zeta$ 
and its is trivial to see that
this means the both  $H(x+iy,a)$ and its first derivative vanishes at $a=0$:

\begin{eqnarray}
H(s,a)_{\vert a=0} =0 , \nonumber\\
\frac{d}{da}H(s,a)_{\vert a=0}=0 .
\label{riemann2}
\end{eqnarray}

\noindent Suppose that  $H(s,a)$  satisfies a differential equation of form

\begin{eqnarray}
\frac{d}{da}H(x+iy,a) = I(x,H(x+iy,a)) ,  
\label{riemann3}
\end{eqnarray}
 
\noindent where $I(x,H)$ is some function having no explicit dependence
on $a$ so that the differential equation defines an autonomous flow.
If the initial conditions  of Eq. \ref{riemann2} are satisfied, 
this differential equation implies that all derivatives of $H$ vanish which in
turn, as it is easy to see,
 implies that the points $s+ m(1-2x)$ are zeros of $\zeta$. This
leaves only the possibility $x=1/2$ so that Riemann hypothesis is proven.
If $I$ is function of also $a$, that is $I=I(a,x,H)$,  this argument breaks down.

The following argument shows that the system is autonomous. One can solve $a$ as function $a=a(x,H)$ from the Taylor series of $H$ with respect to $a$ by using implicit function theorem, substitute this series to the Taylor series of $dH/da$  with respect to $a$,  and by re-organizing the summation 
obtain a Taylor  series with respect to $H$ with coefficients which depend only on $x$ so that one has $I=I(x,H)$.

\subsection{Conclusions}

To sum up, Riemann hypothesis follows from the requirement
that the states in ${\cal{V}}$ can be assigned with a  
conformally invariant physical quantum  system.
This condition reduces to three mutually equivalent conditions: 
the metric induced to ${\cal{V}}$ is
hermitian; positive definite;  allows  conformal  symmetries
as  isometries.   The hermiticity
and positive definiteness properties  reduce to the requirement that the  
dynamical conformal algebra  naturally spanned
by the states in ${\cal{V}}_s$ reduces to the abelian algebra
defined by $L_0=D^+$.  If the  infinitesimal isometries 
for the  matrix elements $\langle \Psi_0\vert \Psi_s\rangle=0$ 
generated by $L_{1-2Re[s]}$ can be exponentiated to isometries as Lie group
theory based  argument strongly suggests, then  Riemann hypothesis follows.

\vm

{\bf Acknowledgements}: I  want  to   express my deep gratitude to Dr. Matthew Watkins for providing me with information about Riemann Zeta and for generous help, in particular for reading the earlier versions of the work and pointing out several inaccuracies and errors. I am also grateful for Prof. Masud Chaichian and   
Doc. Claus Montonen for  encouraging comments and help.

\end{document}